\documentclass[12pt]{article}
\usepackage{amsmath,amssymb,amsbsy}
\usepackage{graphicx}
\usepackage{multirow}
\usepackage{authblk}
\usepackage{algorithm}

\newtheorem{myproposition}{Proposition}[section]
\newtheorem{mytheorem}[myproposition]{Theorem}
\newtheorem{mylemma}[myproposition]{Lemma}

\newtheorem{myobservation}[myproposition]{Observation}

\makeatletter
\def\imod#1{\allowbreak\mkern10mu({\operator@font mod}\,\,#1)}
\makeatother

\begin{document}

\title{Distance magic labeling and two products of graphs\\
}
\author[1]{Marcin Anholcer}
\author[2,$\ast$]{Sylwia Cichacz}
\author[3,$\dagger$]{Iztok Peterin}
\author[3,$\dagger$]{Aleksandra Tepeh}

\let\oldthefootnote\thefootnote
\renewcommand{\thefootnote}{\fnsymbol{footnote}}
\footnotetext[1]{The author was partially supported by National
Science Centre grant nr 2011/01/D/ST/04104, as well as by the Polish
Ministry of Science and Higher Education.} \footnotetext[2]{The
author is also with the Institute  of Mathematics, Physics and
Mechanics, Jadranska 19, 1000 Ljubljana, Slovenia}
\let\thefootnote\oldthefootnote

\affil[1]{\scriptsize{}Pozna\'n University of Economics, Faculty of
Informatics and Electronic Economy} \affil[ ]{Al.
Niepodleg{\l}o\'sci 10, 61-875 Pozna\'n, Poland,
\textit{m.anholcer@ue.poznan.pl}} \affil[ ]{} \affil[2]{AGH
University of Science and Technology, Faculty of Applied
Mathematics} \affil[ ]{Al. Mickiewicza 30, 30-059 Krak\'ow, Poland,
\textit{cichacz@agh.edu.pl}} \affil[ ]{} \affil[3]{University of
Maribor, Faculty of Electrical Engineering and Computer Science}
\affil[ ]{Smetanova 17, 2000 Maribor, Slovenia, \{iztok.peterin,
aleksandra.tepeh\}@uni-mb.si}

\maketitle

\begin{abstract}
Let $G=(V,E)$ be a graph of order $n$. A distance magic labeling of
$G$ is a bijection $\ell \colon V\rightarrow \{1,\ldots ,n\}$ for
which there exists a positive integer $k$ such that $\sum_{x\in
N(v)}\ell (x)=k$ for all $v\in V $, where $N(v)$ is the neighborhood
of $v$. We introduce a natural subclass of distance magic graphs.
For this class we show that it is closed for the direct product with
regular graphs and closed as a second factor for lexicographic
product with regular graphs. In addition, we characterize distance
magic graphs among direct product of two cycles.
\end{abstract}

\affil[1]{\scriptsize{}Pozna\'n University of Economics, Faculty of
Informatics and Electronic Economy} \affil[ ]{Al. Niepodleg{\l}o\'sci 10,
61-875 Pozna\'n, Poland, \textit{m.anholcer@ue.poznan.pl}} \affil[ ]{} %
\affil[2]{AGH University of Science and Technology, Faculty of Applied
Mathematics} \affil[ ]{Al. Mickiewicza 30, 30-059 Krak\'ow, Poland,
\textit{cichacz@agh.edu.pl}} \affil[ ]{} \affil[3]{University of Maribor,
Faculty of Mathematics, Natural Sciences and Information Technologies} %
\affil[ ]{Maribor, Slovenia,
\textit{\{iztok.peterin,aleksandra.tepeh\}@uni-mb.si}} \affil[ ]{} %
\affil[4]{Institute of Mathematics, Physics and Mechanics} \affil[
]{Ljubljana, Slovenia}

\bigskip

\noindent \textbf{Keywords}: distance magic graphs, direct product,
lexicographic product \medskip

\noindent \textbf{AMS subject classification (2010)}: 05C78, 05C76

\section{Introduction and preliminaries}

All graphs considered in this paper are simple finite graphs. We use
$V(G)$ for the vertex set and $E(G)$ for the edge set of a graph
$G$. The \emph{\ neighborhood} $N(x)$ (or more precisely $N_{G}(x)$,
when needed) of a vertex $x$ is the set of vertices adjacent to $x$,
and the \emph{degree} $d(x)$ of $ x$ is $|N(x)|$, i.e. the size of
the neighborhood of $x$. By $C_{n}$ we denote a cycle on $n$
vertices.\newline

\emph{Distance magic labeling} (also called \emph{sigma labeling}) of a
graph $G=(V(G),E(G))$ of order $n$ is a bijection $\ell \colon V\rightarrow
\{1,\ldots ,n\}$ with the property that there is a positive integer $k$
(called \emph{magic constant}) such that $w(x)=\sum_{y\in N_{G}(x)}\ell
(y)=k $ for every $x\in V(G)$, where $w(x)$ is the \emph{weight} of $x $. If
a graph $G$ admits a distance magic labeling, then we say that $G$ is a
\emph{distance magic graph}.\newline

The concept of distance magic  labeling of a graph has been motivated
by the construction of magic squares. However, finding an $r$-regular
distance magic labeling is equivalent to finding equalized
incomplete tournament $\mathrm{EIT}(n, r)$ \cite{FKK1}. In an \emph{equalized incomplete tournament} $\mathrm{EIT}(n, r)$ of $n$ teams with $r$
rounds,  every team plays exactly $r$ other teams and the total
strength of the opponents that team $i$ plays is $k$. For a survey, we refer the reader to \cite{AFK}.\\

The following observations were proved independently:

\begin{myobservation}[\protect\cite{Ji}, \protect\cite{MRS}, \protect\cite%
{Rao}, \protect\cite{Vi}]
\label{obvious}Let $G$ be an $r$-regular distance magic graph on $n$
vertices. Then $k=\frac{r(n+1)}{2}$.
\end{myobservation}

\begin{myobservation}[\protect\cite{Ji}, \protect\cite{MRS}, \protect\cite%
{Rao}, \protect\cite{Vi}]
No $r$-regular graph with an odd $r$ can be a distance magic graph.\label%
{nieparzyste}
\end{myobservation}

We recall three out of four standard graph products (see \cite{IK}). Let $G$
and $H$ be two graphs. All three, the \emph{Cartesian product} $G\square H$,
the \emph{lexicographic product} $G\circ H$, and the \emph{direct product} $
G\times H$ are graphs with vertex set $V(G)\times V(H)$. Two vertices $(g,h)$
and $(g^{\prime },h^{\prime })$ are adjacent in:

\begin{itemize}
\item $G\square H$ if and only if either $g=g^{\prime }$ and $h$ is adjacent
with $h^{\prime }$ in $H$, or $h=h^{\prime }$ and $g$ is adjacent with $
g^{\prime }$ in $G$;

\item $G\circ H$ if and only if either $g$ is adjacent with $g^{\prime }$ in
$G$ or $g=g^{\prime }$ and $h$ is adjacent with $h^{\prime }$ in $H$;

\item $G\times H$ if and only if $g$ is adjacent with $g^{\prime }$ in $G$
and $h$ is adjacent with $h^{\prime }$ in $H$.
\end{itemize}

For a fixed vertex $g$ of $G$, the subgraph of any of the above products
induced by the set $\{(g,h)\,:\,h \in V(H)\}$ is called an \emph{$H$-layer}
and is denoted $^g\!H$. Similarly, if $h \in H$ is fixed, then $G^h$, the
subgraph induced by $\{(g,h)\,:\, g \in V(G)\}$, is a \emph{$G$-layer}.

The main topic of this paper is the direct product (that is known also by
many other names, see \cite{IK}). It is the most natural graph product in
the sense that each edge of $G\times H$ projects to an edge in both factors $
G$ and $H$. This is also the reason that many times this product is the most
difficult to handle among (standard) products. Even the distance formula is
very complicated with respect to other products (see \cite{Kim}) and $%
G\times H $ does not need to be connected, even if both factors are.
More precisely, $G\times H$ is connected if and only if both $G$ and
$H$ are connected and at least one of them is non-bipartite
\cite{Weich}.

The direct product is commutative, associative, and has attracted a
lot of attention in the research community in last 50 years.
Probably the biggest challenge (among
all products) is the famous Hedetniemi's conjecture:%
\begin{equation*}
\chi (G\times H)=\min \{\chi (G),\chi (H)\}.
\end{equation*}%
This conjecture suggests that the chromatic number of the direct product
depends only on the properties of one factor and not both. This is not so
rare and also in this work we show that it is enough for one factor to be a
distance magic graph with one additional property and then the product with
any regular graph will result in a distance magic graph. For more about the
direct product and products in general we recommend the book \cite{IK}.

Some graphs which are distance magic among
(some) products can be seen in \cite{Be,Cic,MRS,RSP}. The following product cycle and product related results were proved by
Miller, Rodger, and Simanjuntak.

\begin{mytheorem}[\protect\cite{MRS}]
\label{MRSC4} The cycle $C_n$ of length $n$ is a distance magic graph if and
only if $n = 4 $.
\end{mytheorem}

\begin{mytheorem}[\protect\cite{MRS}]
\label{lex_pr} Let $G$ be an $r$-regular graph and $C_{n}$ the cycle of
length $n$ for $r\geq 1$ and $n\geq 3$. The lexicographic product $G\circ
C_{n}$ admits a distance magic labeling if and only if $n=4$.
\end{mytheorem}

In particular we have:

\begin{myobservation}
The lexicographic product $C_n \circ C_m$, $n, m\geq3$ is a distance magic
graph if and only if $m =4$.
\end{myobservation}

Rao, Singh and Parameswaran characterized distance magic graphs
among Cartesian products of cycles.

\begin{mytheorem}[\protect\cite{RSP}]
\label{cart_cycle1}The Cartesian product $C_{n}\square C_{m}$, $n,m\geq 3$,
is a distance magic graph if and only if $n=m\equiv 2\imod 4)$.
\end{mytheorem}

In the next section we introduce a natural subclass of distance magic
graphs. For this class of graphs we were able to generalize the Theorem \ref%
{lex_pr} and show that it is closed for the direct product with
regular graphs. In the last section we characterize distance magic
graphs among direct products of cycles. In particular, we prove that
a graph $C_{m}\times C_{n}$ is distance magic if and only if $n=4$
or $m=4$ or $m,n\equiv 0 \imod 4) $.


\section{Balanced distance magic graphs}

In order to obtain a large class of graphs for which their direct product is
a distance magic graph we introduce a natural subclass of distance magic
graphs.

A distance magic graph $G$ with an even number of vertices is
called\emph{\ balanced} if there exists a bijection $\ell
:V(G)\rightarrow \{1,\ldots ,|V(G)|\}$ such that for every $w\in
V(G)$ the following holds: if $u\in N(w) $ with $\ell (u)=i$, there
exists $v\in N(w)$ with $\ell (v)=|V(G)|+1-i$ . We call $u$ the
\emph{twin vertex} of $v$ and vice versa (we will also say that $u$
and $v$ are \emph{twin vertices}, or shortly \emph{twins}) and $\ell
$ is called a \emph{\ balanced distance labeling}. Hence a distance
magic graph $G$ is balanced if for any $w\in V(G)$ either both or
none of vertices $u$ and $v$ with labels $ \ell (u)=i$ and $\ell
(v)=|V(G)|+1-i$ are in the neighborhood of $w$ . It also follows
from the definition that twin vertices of a balanced distance magic
graph cannot be adjacent and that $N_{G}(u)=N_{G}(v)$.

It is somewhat surprising that the condition $N_{G}(u)=N_{G}(v)$
plays an important role in finding the factorization of the direct
product, see Chapter 8 of \cite{IK}. In particular, if a
non-bipartite connected graph has no pairs of vertices with the
property $N_G(u)=N_G(v)$, then it is easier to find the prime factor
decomposition. Similarly, such pairs generate very simple
automorphisms of $G$ and have been called unworthy in \cite{Wils}.
However in both above mentioned cases not all vertices need to have
a twin vertex as in our case.

It is easy to see that a balanced distance magic graph is an $r$-regular
graph for some even $r$. Recall that the magic constant is $\frac{r}{2}%
(|V(G)|+1)$ by Observation \ref{obvious}. Trivial examples of balanced
distance magic graph are graphs with no edges and even number of vertices.
Not all distance magic graphs are balanced distance magic graphs. The
smallest example is $P_{3}$. More examples (regular graphs with an even
number of vertices) will be presented in next section.

The graph $K_{2n,2n}$, $n\geq 1$, is a balanced distance magic
graph. To verify this let $V(K_{2n,2n})=\{v_{1},\ldots ,v_{4n}\}$.
Assume that the vertices are enumerated in such a way that the sets
$U=\{v_{i}:i(\mathop{\rm mod}\nolimits 4)\in \{0,1\}\}$ and
$W=V(K_{2n,2n})-U$ form the bipartition of $V(K_{2n,2n})$. It is
easy to see that the labeling

\begin{equation*}
\ell(v_{i})=i \text{ for } i\in \{1,\ldots ,4n\}
\end{equation*}

\noindent{}is the desired balanced distance magic labeling for $n\geq 2$. In
particular, for $n=1$ note that $K_{2,2}$ is isomorphic to $C_{4}$ and
consecutive vertices receive labels $1,2,4,3$.

Also $K_{2n}-M$ is a balanced distance magic graph if $M$ is a perfect
matching of $K_{2n}$. Indeed, if $u$ and $v$ form an $i$-th edge of $M$, $
i\in \{1,\ldots ,n\}$, we set $\ell (u)=i$ and $\ell (v)=2n+1-i$ which is a
balanced distance magic labeling.

The distance magic graphs $G\circ C_{4}$ described in Theorem \ref{lex_pr}
are also balanced distance magic graphs. Let $V(G)=\{g_{1},\ldots ,g_{p}\}$
be the vertex set of a regular graph $G$ and $V(C_{4})=
\{h_{1},h_{2},h_{3},h_{4}\}$ where indices of vertices in $V(C_{4})$
correspond to labels of a distance magic labeling of $C_{4}$. It is not hard
to verify that the labeling
\begin{equation*}
\ell ((g_{i},h_{j}))=\left\{
\begin{array}{lcl}
(j-1)p+i, & \text{if} & j\in \{1,2\}, \\
jp-i+1, & \text{if} & j\in \{3,4\},%
\end{array}%
\right.
\end{equation*}%
is a balanced distance magic labeling of $G\circ C_{4}$. Using similar
labeling we obtain a larger family of balanced distance magic graphs.


\begin{mytheorem}
\label{lex-BDM} Let $G$ be a regular graph and $H$ a graph not isomorphic to
$\overline{K_{n}}$ where $n$ is odd. Then $G\circ H$ is a balanced distance
magic graph if and only if $H$ is a balanced distance magic graph.
\end{mytheorem}

\noindent \textit{Proof.}\ Let $G$ be an $r_G$-regular graph and $H$
a graph not isomorphic to $\overline{K_{n}}$ for an odd $n$. Let
first $H$ be a balanced distance magic graph with the vertex set
$V(H)=\{h_{1},\ldots ,h_{t}\}$ and let $\varphi $ defined by
$\varphi (h_{j})=j$ be a balanced distance magic labeling of $H$ (we
can always enumerate the vertices in an appropriate way). Recall
that $t=|V(H)|$ is an even number, $H$ is an $r_{H}$ -regular graph
where $r_{H}$ is also even and that for $j\leq \frac{t}{2}$,
$h_{t+1-j} $ is the twin vertex of $h_{j}$. Let $V(G)=\{g_{1},\ldots
,g_{p}\}$.

For $i\in \{1,\ldots ,p\}$ and $j\in \{1,\ldots ,t\}$ define the following
labeling $\ell $:{}

\begin{equation*}
\ell (g_{i},h_{j})=\left\{
\begin{array}{lcl}
(j-1)p+i, & \text{if} & j\leq \frac{t}{2}, \\
jp-i+1, & \text{if} & j>\frac{t}{2}.%
\end{array}
\right.
\end{equation*}

It is straightforward to see that $\ell$ is a bijection. For $j\leq \frac{t%
}{2}$ we have
$\ell(g_i,h_j)+\ell(g_i,h_{t+1-j})=(j-1)p+i+(t+1-j)p-i+1=tp+1$. The
fact that $H$ is a balanced distance magic graph and the structure
of the graph $G\circ H$ together imply that if $(g_i,h_j)$ is a
neighbor of some vertex $(g,h)\in V(G\circ H)$ then also
$(g_i,h_{t+1-j})$ is a neighbor of this vertex. We derive that
$(g_i,h_j)$ and $(g_i,h_{t+1-j})$ are twin vertices.

To finish the proof that $G\circ H$ is a balanced distance magic graph we
now only need to verify that the weights of all the vertices $(g,h)$ in $%
G\circ H$ are equal:

\begin{eqnarray*}
w(g,h) &=&\sum_{(g_{i},h_{j})\in N_{G\circ H}((g,h))}\ell (g_{i},h_{j})= \\
&=&\sum_{g_{i}\in N_{G}(g)}\sum_{h_{j}\in V(H)}\ell
(g_{i},h_{j})+\sum_{h_{j}\in N_{H}(h)}\ell (g,h_{j}) = \\
&=&r_{G}\frac{t}{2}(tp+1)+\frac{r_{H}}{2}(tp+1)=\frac{(tr_{G}+r_{H})(tp+1)}{2%
}.
\end{eqnarray*}

Conversely, let $G\circ H$ be a balanced distance magic (and hence a
regular) graph. If $H$ is an empty graph on even number of vertices, it is
balanced distance magic graph.

In the case when $H$ is not an empty graph we claim that the twin vertex of
any $(g,h)\in V(G\circ H)$ lies in $^g\!H$. Suppose to the contrary that
there exist twin vertices $(g,h)$ and $(g^{\prime},h^{\prime})$ such that $%
g\neq g^{\prime}$. Then $g$ and $g^{\prime}$ are at distance $2$ in
$G$ ($gg^{\prime}\in E(G)$ would imply that $(g,h)$ and
$(g^{\prime},h^{\prime})$ are adjacent, which is impossible).
Assumption that $h$ is not an isolated vertex in $H$ leads to a
contradiction, since if there is an edge $hh^{\prime\prime}\in
E(H)$,
then $(g,h^{\prime\prime})\in N_{G\circ H}((g,h))$ but $(g,h^{\prime\prime})%
\notin N_{G\circ H}((g^{\prime},h^{\prime}))$ (recall that twin vertices
have the same neighborhood). Since $H$ is a regular graph (it is easy tho
see that if it was not, then $G\circ H$ would not be regular either) we
derive that $H$ is an empty graph, a contradiction. Thus two twin vertices
of $G\circ H$ lie in the same $H$-layer.

This implies that $H$ has an even number of vertices $t=|V(H)|$. Let
$V(H)=\{v_{1},\ldots ,v_{\frac{t}{2}},v_{1}^{\prime
},\ldots ,v_{\frac{t}{2}}^{\prime }\}$ where we use this notation to indicate that $(g,v_{i})$ and $%
(g,v_{i}^{\prime })$ are the twin vertices in $^{g}\!H$, $i\in \{1,\ldots ,%
\frac{t}{2}\}$. To prove that $H$ is a balanced distance magic graph we need
to see that the function $\ell :V(H)\rightarrow \{1,\ldots ,t\}$ defined by $%
\ell (v_{i})=i$ and $\ell (v_{i}^{\prime })=t-i+1$ for $i\in \{1,\ldots ,%
\frac{t}{2}\}$ is a magic distance labeling of $H$.

Obviously, $\ell $ is a bijection. As $H$ is a regular, nonempty graph, each
pair of twin vertices $(g,v_{i})$ and $(g,v_{i}^{\prime })$ appears in the
neighborhood of some vertex $(g,u)$, where $u\neq v_{i}$ and $u\neq
v_{i}^{\prime }$, thus $v_{i}$ and $v_{i}^{\prime }$ are neighbors of $u$ in
$H$. Since $H$ is an $r_{H}$-regular graph, we deduce that every vertex $v$
in $H$ has $r_H/2$ pairs $(v_{i},v_{i}^{\prime })$ of neighbors, and each
such pair contributes $t+1$ to the weight of $v$. Hence $w(v)=\frac{%
r_{H}(t+1)}{2}$ and $H$ is a balanced distance magic graph.\hfill \rule%
{0.1in}{0.1in}\medskip

Note that in order to prove the equivalence in the above theorem we
needed to exclude $H$ as an empty graph with odd number of vertices.
Namely, it
is not hard to see that for positive integer $k$, $C_4\circ \overline{%
K_{2k-1}}$ is a balanced distance magic graph, but
$\overline{K_{2k-1}}$ is not (recall that by the definition an empty
graph is balanced distance magic if it has an even order). As an
example see the labeling of $C_4\circ \overline{K_{3}}$ in the table
below, where rows and columns represent labeling of vertices in
$C_{4}$-layers and $\overline{K_{3}}$-layers, respectively (the
latter ones refer to consecutive vertices of $C_4$).

\begin{tabular}[t]{|c|c|c|c|}
\multicolumn{4}{c}{} \\ \hline
3 & 6 & 10 & 7 \\ \hline
2 & 5 & 11 & 8 \\ \hline
1 & 4 & 12 & 9 \\ \hline
\end{tabular}%
\medskip

The situation is even more challenging when we turn to the direct
product. If one factor, say $H$, is an empty graph, also the product
$G\times H$ is an empty graph. Hence for any graph $G$ on even
number of vertices $G\times
\overline{K_{2k-1}}$ is a balanced distance magic graph, while $\overline{%
K_{2k-1}}$ is not. However, we can still obtain the result only
slightly weaker than Theorem \ref{lex-BDM}. For this we need the
following observations.

\begin{mylemma}
\label{lemma1} Let $G\times H$ be a balanced distance magic graph
and let $(g,h)$ and $(g',h')$ with $g\neq g'$ and $h\neq h'$ be twin
vertices for some balanced distance magic labeling. The labeling in
which we exchange the labels of $(g',h')$ and $(g',h)$ is a balanced
distance magic labeling with $(g,h)$ and $(g',h)$ as twin vertices.
\end{mylemma}

\noindent \textit{Proof.}\ Let $\ell: V(G\times H)\rightarrow
\{1,\ldots,|V(G)||V(H)|\}$ be a balanced distance magic labeling
where $(g,h)$ and $(g',h')$ are twin vertices with $g\neq g'$ and
$h\neq h'$. Recall that $N_{G \times H}(a,b)=N_G(a)\times N_H(b)$
for every $(a,b)\in V(G \times H)$ and that twin vertices have the
same neighborhood. Thus we derive
\begin{equation*}
N_{G\times H}(g,h)=N_{G\times H}(g',h')=N_{G\times H}(g',h)=N_{G\times H}(g,h').
\end{equation*}

Using this property we can show that the labeling $ \widehat{\ell}:
V(G\times H) \rightarrow \{1,\ldots,|V(G)||V(H)|\}$ defined by
$\widehat{\ell}(g',h)=\ell(g',h')$,
$\widehat{\ell}(g',h')=\ell(g',h)$ and
$\widehat{\ell}(a,b)=\ell(a,b)$ for every $(a,b)\in V(G\times
H)\setminus \{(g',h'), (g',h)\}$ is a balanced distance magic
labeling of $V(G\times H)$. To show this let $(g'',h'')$ be the twin
vertex of $(g',h)$, and $(g''',h''')$ the twin vertex of $(g,h')$
with respect to the labeling $\ell$.

If $(a,b)$ is a vertex that is not adjacent to any vertex in
$$S=\{(g,h),(g,h'),(g',h),(g',h'),(g'',h''),(g''',h''')\},$$ then the
label under $\widehat{\ell}$ of every neighbor of $(a,b)$ remains
the same as under $\ell$ and since $\ell$ is a balanced distance
magic labeling every vertex in $N_{G\times H}(a,b)$ has its twin
vertex in $N_{G\times H}(a,b)$.
(Note that also the case when $(a,b)\in S$ is
included here.)

If $(a,b)$ is adjacent to at least one vertex from $S$, one can
observe that $(a,b)$  is in fact adjacent to all vertices in $S$.
Hence also in this case we derive that every
vertex in the open neighborhood of $(a,b)$ has its twin vertex
within this neighborhood.

Since $\widehat{\ell}$ is obviously a bijection we have proved that
$\widehat{\ell}$ is a balanced distance magic labeling where $(g,h)$
and $(g',h)$ are twin vertices. \hfill \rule{0.1in}{0.1in}\medskip

This lemma has clearly a symmetric version if we exchange the labels of
$(g',h')$ and $(g,h')$.

\begin{mylemma}
\label{lemma2} Let $G\times H$ be a balanced distance magic graph,
and let $(g,h)$ and $(g',h)$ be twin vertices as well as $(g,h_1)$
and $(g,h_2)$ for some balanced distance magic labeling. The
labeling in which we exchange the labels of $(g,h_2)$ and $(g',h_1)$
is a balanced distance magic labeling with twins $(g,h_1)$ and
$(g',h_1)$.
\end{mylemma}

\noindent \textit{Proof.}\ Let
$\ell: V(G\times H)\rightarrow \{1,\ldots,|V(G)||V(H)|\}$ be a balanced distance
magic labeling of $G\times H$ where $\{(g,h),(g',h)\}$ and $\{(g,h_1),(g,h_2)\}$
are pairs of twin vertices. As in the proof of Lemma \ref{lemma1} we have
\begin{equation*}
N_{G\times H}(g,h_1)=N_{G\times H}(g,h_2)=N_{G\times H}(g',h_1)=N_{G\times H}(g',h_2).
\end{equation*}
By the same arguments as in the proof of Lemma \ref{lemma1} it is
easy to see that the labeling  $ \widehat{\ell}:V(G\times H)
\rightarrow \{1,\ldots,|V(G)||V(H)|\}$ defined by
$\widehat{\ell}(g,h_2)=\ell(g',h_1)$,
$\widehat{\ell}(g',h_1)=\ell(g,h_2)$ and
$\widehat{\ell}(a,b)=\ell(a,b)$ for every $(a,b)\in V(G\times
H)\setminus \{(g',h_1), (g,h_2)\}$ is a balanced distance magic
labeling of $V(G\times H)$. Clearly, $(g,h_1)$ and $(g',h_1)$ are
twins for $ \widehat{\ell}$. \hfill \rule{0.1in}{0.1in}\medskip

\begin{mylemma}
\label{lemma3} Let $G\times H$ be a balanced distance magic graph,
and let $(g,h)$ and $(g',h)$ be twin vertices as well as $(g,h')$
and $(g'',h')$, $g''\neq g'$, for some balanced distance magic
labeling. The labeling in which we exchange the labels of $(g',h')$
and $(g'',h')$ is a balanced distance magic labeling where $(g,h')$
and $(g',h')$ are twin vertices.
\end{mylemma}

\noindent \textit{Proof.}\ Let $\ell: V(G\times H)\rightarrow
\{1,\ldots,|V(G)||V(H)|\}$ be a balanced distance magic labeling of
$G\times H$ where $\{(g,h),(g',h)\}$ and $\{(g,h'),(g'',h')\}$ are
pairs of twin vertices for $g''\neq g'$. One can observe that
\begin{equation*}
N_{G\times H}(g,h')=N_{G\times H}(g',h')=N_{G\times H}(g'',h').
\end{equation*}
Using similar arguments as in the proof of Lemma \ref{lemma1} it is
easy to see that the labeling  $ \widehat{\ell}:V(G\times H)
\rightarrow \{1,\ldots,|V(G)||V(H)|\}$ defined by
$\widehat{\ell}(g',h')=\ell(g'',h')$,
$\widehat{\ell}(g'',h')=\ell(g',h')$ and
$\widehat{\ell}(a,b)=\ell(a,b)$ for every $(a,b)\in V(G\times
H)\setminus \{(g',h'), (g'',h')\}$ is a balanced distance magic
labeling of $V(G\times H)$. Clearly, $(g,h')$ and $(g',h')$ are
twins with respect to the labeling $ \widehat{\ell}$. \hfill
\rule{0.1in}{0.1in}\medskip

\begin{mytheorem}
\label{str_pr} The direct product $G\times H$ is a balanced distance
magic graph if and only if one of the graphs $G$ and $H$ is a
balanced distance magic and the other a regular graph.
\end{mytheorem}

\noindent \textit{Proof.}\ Assume first, without loss of generality,
that G is a regular and $H$ is a balanced distance magic graph with
$V(H)=\{h_{1},\ldots ,h_{p}\}$, where a suffix
indicates the label of a balanced distance magic labeling of $H$. Thus for $%
j\leq \frac{p}{2}$, $h_{p+1-j}$ is the twin vertex of $h_{j}$. Recall that $r_{H}$ is even. Let $V(G)=\{g_{1},\ldots
,g_{t}\}$.

For $i\in \{1,\ldots ,t\}$ and $j\in \{1,\ldots ,p\}$ define the following
labeling $\ell $:{}

\begin{equation*}
\ell (g_{i},h_{j})=\left\{
\begin{array}{lcl}
(j-1)t+i, & if & j\leq \frac{p}{2}, \\
jt-i+1, & if & j>\frac{p}{2}.%
\end{array}
\right.
\end{equation*}

\noindent{}It is straightforward to see that
$\ell :V(G\times H)\rightarrow \{1,\ldots ,pt\}$ is a bijection.
Moreover, note that for any $j\leq \frac{p}{2}$ we have $\ell
(g_{i},h_{j})+\ell (g_{i},h_{p+1-j})=(j-1)t+i+(p+1-j)t-i+1=pt+1$. Moreover, if $%
(g_{i},h_{j})\in N_{G\times H}(g,h)$, then also
$(g_{i},h_{p+1-j})\in N_{G\times H}(g,h)$, since $h_{j}\in N_{H}(h)$
implies that $h_{p+1-j}\in N_{H}(h)$. Hence $ (g_{i},h_{j})$ is the
twin vertex of $(g_{i},h_{p+1-j})$.

Finally, we finish the proof of the first implication by the
following calculation for an arbitrary vertex $(g,h)\in V(G\times H)$:%
\begin{eqnarray*}
w(g,h) &=&\sum_{(g_{i},h_{j})\in N_{G}(g)\times N_{H}(h)}\ell
(g_{i},h_{j})=\sum_{g_{i}\in N_{G}(g)}\sum_{h_{j}\in N_{H}(h)}\ell
(g_{i},h_{j})= \\
&=&\sum_{g_{i}\in N_{G}(g)}\sum_{h_{j}\in N_{H}(h),\ j\leq \frac{p}{2}}(\ell
(g_{i},h_{j})+\ell(g_{i},h_{p+1-j})) = \\
&=&\sum_{g_{i}\in N_{G}(g)}\sum_{h_{j}\in N_{H}(h),\ j\leq \frac{p}{2}%
}(pt+1)= \\
&=&\frac{r_{H}r_{G}}{2}(pt+1).
\end{eqnarray*}

Conversely, let $G\times H$ be a balanced distance magic graph (this
implies that $G\times H$ is a regular graph and hence also both $G$
and $H$ are regular). There exists a balanced distance magic
labeling $\ell: V(G\times H) \rightarrow \{1,\ldots,|V(G)||V(H)|\}$.
First we show the following.

\noindent\textbf{Claim}
\textit{There exists a balanced distance labeling of $G\times H$ such that
one of the following is true:}
\begin{enumerate}
\item
\textit{There exists an $H$-layer $^{g}\!H$, such that the twin vertex of any $(g,h)\in {^{g}\!H}$ lies in $^{g}\!H$.}
\item
\textit{There exists a $G$-layer $G^{h}$, such that the twin vertex of any $(g,h)\in {G^{h}}$ lies in $G^{h}$.}
\end{enumerate}

If there exists an $H$-layer or a $G$-layer such that the twin
vertex of any vertex in this layer also lies within this layer, then
we are done. Hence assume that this is not the case, i.e. for every
$H$-layer $^{g}\!H$ there exists a vertex $(g,h)$ such that its twin
vertex $(g',h')$ has the property $g'\neq g$.

We use Algorithm \ref{alg_couple} to rearrange the labels of
vertices in $V(G\times H)$ in such a way that we either obtain an
$H$-layer closed for twin vertices or we couple all the $H$-layers,
i.e. we find pairs of $H$-layers $\{^{g}\!H,{^{g'}\!H}\}$ with the
property that the twin vertex of a vertex $(g,h)\in ^g\!H$ lies in
${^{g'}\!H}$ and is of the form $(g',h)$. The latter case implies
that all the $G$-layers (and in particular one of them, say $G^h$)
are closed for twins, and the claim is proved.

\begin{algorithm}
\caption{Coupling $H$-layers}\label{alg_couple}
\begin{enumerate}
\item[Step 1:]
Set $A=V(G)$. Go to step $2$.
\item[Step 2:]
If $A=\{g\}$ for some $g$, then STOP, $^{g}\!H$ is closed under twin
vertices. If $A=\emptyset$ then STOP, all the $H$-layers are matched
in such a way that in any pair $\{^{g}\!H,{^{g'}\!H}\}$ for every
vertex $(g,h)$ its twin vertex is of the form $(g',h)$. If $|A|\geq
2$, then proceed to step 3.
\item[Step 3:]
Choose any $g\in A$. If $^{g}\!H$ is closed for twin vertices, then
STOP. Otherwise, there is a vertex $(g,h)\in {^{g}\!H}$ having the
twin $(g',h')$, where $g'\in A$ and $g'\neq g$. If $h'\neq h$, then
use Lemma \ref{lemma1} to obtain a new labeling with $(g,h)$ and
$(g',h)$ being twins. If every vertex $(a,b)\in {^{g}\!H}\cup
{^{g'}\!H}$ has its twin vertex $(a',b')$ also in ${^{g}\!H}\cup
{^{g'}\!H}$, then go to step 6. Otherwise go to step 4.
\item[Step 4:]
For every vertex $(g,h_1)\in {^{g}\!H}$ with the twin vertex
$(g'',h_2)$, where $g''\notin\{g,g'\}$, $h_2\neq h_1$, use Lemma
\ref{lemma1} to obtain a new labeling where $(g,h_1)$ and
$(g'',h_1)$ are twin vertices. Go to step 5.
\item[Step 5:]
For every vertex $(g,h_1)\in {^{g}\!H}$ with the twin vertex
$(g'',h_1)$, where $g''\notin\{g,g'\}$, use Lemma \ref{lemma3} to
obtain a new labeling with twin vertices $(g,h_1)$ and $(g',h_1)$.
Proceed to step 6.
\item[Step 6:]
Until there exists a pair of twin vertices $(g,h_1),(g,h_2)\in
{^{g}\!H}$ with the property $h_2\neq h_1$, use Lemma \ref{lemma2}
to obtain a new labeling where $(g,h_1)$ and $(g',h_1)$ are twin
vertices. Proceed to step 7.
\item[Step 7:]
Set $A\gets A\setminus\{g,g'\}$ and go back to step 2.
\end{enumerate}
\end{algorithm}

Assume that some $H$-layer, say $^{g}\!H$, is closed for twins. In
this case $H$ has an even number of vertices $p=|V(H)|$. We
enumerate the vertices as follows $V(H)=\{v_{1},\ldots
,v_{\frac{p}{2}},v_{1}^{\prime },\ldots ,v_{\frac{p}{2}}^{\prime
}\}$ in such a way that $(g,v_{i})$ and $(g,v_{i}^{\prime })$ are
twin vertices in $^{g}H$, for $i\in \{1,\ldots ,\frac{p}{2}\}$. To
prove that $H$ is a balanced distance magic graph we need to see
that the function $\ell :V(H)\rightarrow \{1,\ldots ,p\}$ defined by
$\ell (v_{i})=i$ and $\ell (v_{i}^{\prime })=p-i+1$ for $i\in
\{1,\ldots ,\frac{p}{2}\}$ is a magic distance labeling of $H$.

Obviously, $\ell $ is a bijection. Note that any pair of twin vertices $%
(g,v_{i})$ and $(g,v_{i}^{\prime })$ appears in the neighborhood of some
vertex $(g,u)$, where $u\neq v_{i}$ and $u\neq v_{i}^{\prime }$, thus $v_{i}$
and $v_{i}^{\prime }$ are both neighbors of $u$ in $H$. Since $H$ is an $%
r_{H}$-regular graph, we deduce that every vertex $v$ in $H$ has $r_H$ pairs
$(v_{i},v_{i}^{\prime })$ of neighbors, and each such pair contributes $p+1$
to the weight of $v$. Hence $w(v)=\frac{r_{H}(p+1)}{2} $ and $H$ is a
balanced distance magic graph.

In the case when some $G$-layer $G^{h}$ is closed for twins, we can
prove in an analogous way that $G$ is a balanced distance magic
graph. ~\hfill \rule{0.1in}{0.1in}\medskip


\section{Distance magic graphs $C_{m}\times C_{n}$}

Let $V(C_{m}\times C_{n})=\{v_{i,j}:0\leq i\leq m-1,0\leq j\leq
n-1\}$, where
$N(v_{i,j})=\{v_{i-1,j-1},v_{i-1,j+1},v_{i+1,j-1},v_{i+1,j+1}\}$ and
operation on the first suffix is taken modulo $m$ and on the second
suffix modulo $n$. We also refer to the set of all vertices
$v_{i,j}$ with fixed $i$ as $i$-th row and with fixed $j$ as $j$-th
column.

We start with direct products of cycles that are not distance magic.

\begin{mytheorem}
\label{notdm}If $m\not\equiv 0\imod 4$ and $n\neq 4$ or $n\not\equiv
0\imod 4$ and $m\neq 4$, then $C_{m}\times C_{n}$ is not distance
magic.
\end{mytheorem}

\noindent \textit{Proof.}\ By commutativity of the direct product we
can assume that $m\not\equiv 0\imod 4$ and $n\neq 4$. Assume that
$C_{m}\times C_{n}$ is distance magic with some magic constant $k$,
which means there is a distance magic labeling $\ell$. Let us
consider the neighborhood sum of labels of $v_{i+1,j+1}$ and
$v_{i+3,j+1}$ for any $i\in\{0,\dots,m-1\}$ and
$j\in\{0,\dots,n-1\}$:
\begin{equation*}
w(v_{i+1,j+1})=\ell (v_{i,j})+\ell (v_{i,j+2})+\ell (v_{i+2,j})+\ell
(v_{i+2,j+2})=k,
\end{equation*}
\begin{equation*}
w(v_{i+3,j+1})=\ell (v_{i+2,j})+\ell (v_{i+2,j+2})+\ell (v_{i+4,j})+\ell
(v_{i+4,j+2})=k.
\end{equation*}
It implies that
\begin{equation*}
\ell (v_{i,j})+\ell (v_{i,j+2})=\ell (v_{i+4,j})+\ell (v_{i+4,j+2}).
\end{equation*}
Repeating that procedure we obtain that
\begin{equation*}
\ell (v_{i,j})+\ell (v_{i,j+2})=\ell (v_{i+4\alpha ,j})+\ell (v_{i+4\alpha
,j+2})
\end{equation*}%
for any natural number $\alpha $.

It is well known that if $a,b\in \mathop\mathbb{Z}\nolimits_m$ and $%
\mathop{\rm gcd}\nolimits(a,m)=\mathop{\rm gcd}\nolimits(b,m)$, then $a$ and
$b$ generate the same subgroup of $\mathop\mathbb{Z}\nolimits_m$, that is, $%
\langle a\rangle=\langle b\rangle$.

Since $m\not\equiv 0\imod 4$ we
have $\mathop{\rm gcd} \nolimits(2,m)=\mathop{\rm gcd}\nolimits(4,m)$ and $%
2\in \langle 4\rangle $, which implies that there exists $\alpha
^{\prime }$ such that $4\alpha ^{\prime }\equiv 2\imod m$. We deduce
that
\begin{equation*}
\ell(v_{i,j})+\ell (v_{i,j+2})=\ell (v_{i+2,j})+\ell (v_{i+2,j+2})=\frac{k}{2%
}.
\end{equation*}

Substituting $j$ with $j+2$ we obtain
\begin{equation*}
\ell (v_{i,j+2})+\ell (v_{i,j+4})=\frac{k}{2}.
\end{equation*}

\noindent{}Thus for every $i, j$ we have
\begin{equation*}
\ell(v_{i,j})=\ell(v_{i,j+4}),
\end{equation*}
which leads to a contradiction, since $n\neq 4$ and $\ell$ is not a
bijection.~\hfill \rule{0.1in}{0.1in}\medskip

Next we show that some of direct products of cycles are distance
magic but not balanced distance magic. Used constructions are
similar to those by Cichacz and Froncek in \cite{CicFro}.

\begin{mytheorem}
\label{0mod4}If $m,n\equiv 0\imod 4$
and $m,n>4$, then $C_{m}\times C_{n}$ is distance magic but not balanced
distance magic graph.
\end{mytheorem}

\noindent \textit{Proof.}\ First we show that $C_{m}\times C_{n}$ is
distance magic. We define the labeling $\ell $ by starting conditions (every
second vertex of the row zero) followed by recursive rules that cover all
the remaining vertices.

\begin{equation*}
\ell (v_{0,4j+t})=\left\{
\begin{array}{lclcl}
2j+1, & \text{if} & 0\leq j \leq \lceil\frac{n}{8}\rceil-1 & \text{and} &
t=0, \\
\frac{n}{2}-2j, & \text{if} & \lceil\frac{n}{8}\rceil\leq j \leq \frac{n}{4}%
-1 & \text{and} & t=0, \\
mn-2j-1, & \text{if} & 0\leq j \leq \lfloor\frac{n}{8}\rfloor-1 & \text{and}
& t=2, \\
mn-\frac{n}{2}+2j+2, & \text{if} & \lfloor\frac{n}{8}\rfloor\leq j \leq
\frac{n}{4}-1 & \text{and} & t=2.%
\end{array}
\right.
\end{equation*}

Note that we have used every label between $1$ and $\frac{n}{4}$ as well as
between $mn-\frac{n}{4}+1$ and $mn$ exactly once for the starting conditions.

In the first recursive step we label every second vertex of row two
in the order that is in a sense opposite to the one of row zero:

\begin{equation*}
\ell (v_{2,j})=\left\{
\begin{array}{lcl}
\ell (v_{0,n-2-j})+\frac{n}{4}, & \text{if} & \ell (v_{0,n-2-j})\leq \frac{mn%
}{2}, \\
\ell (v_{0,n-2-j})-\frac{n}{4}, & \text{if} & \ell (v_{0,n-2-j})>\frac{mn}{2}%
,%
\end{array}
\right.
\end{equation*}
for $j\in \{0,2,\ldots ,n-2\}$. Clearly we use in this step every label
between $\frac{n}{4}+1$ and $\frac{n}{2}$ and between $mn-\frac{n}{2}+1$ and
$mn-\frac{n}{4}$ exactly once.

We continue with every second vertex in every even row. Hence for $2\leq
i\leq \frac{m}{2}-1$ and for $j\in \{0,2,\ldots ,n-2\}$ let

\begin{equation*}
\ell (v_{2i,j})=\left\{
\begin{array}{lcl}
\ell (v_{2i-4,j})+\frac{n}{2}, & \text{if} & \ell (v_{2i-4,j})\leq \frac{mn}{%
2}, \\
\ell (v_{2i-4,j})-\frac{n}{2}, & \text{if} & \ell (v_{2i-4,j})>\frac{mn}{2}.%
\end{array}
\right.
\end{equation*}%
Again all labels here are used exactly once and are between $\frac{n}{2}+1$
and $\frac{mn}{8}$ and between $mn-\frac{mn}{8}+1=\frac{7mn}{8}+1$ and $mn-%
\frac{n}{2}$.

Next we label every second vertex of every odd row and complete with this
all even columns. For $0\leq i\leq \frac{m}{2}-1$ and for $j\in
\{0,2,\ldots,n-2\}$ we set:

\begin{equation*}
\ell (v_{2i+1,j})=\left\{
\begin{array}{lcl}
\ell (v_{2i,j})+\frac{mn}{8} & \text{if} & \ell (v_{2i,j})\leq \frac{mn}{2},
\\
\ell (v_{2i,j})-\frac{mn}{8} & \text{if} & \ell (v_{2i,j})>\frac{mn}{2}.%
\end{array}
\right.
\end{equation*}
Labels used here are between $\frac{mn}{8}+1$ and $\frac{mn}{4}$ and between
$\frac{3mn}{4}+1$ and $\frac{7mn}{8}$.

Finally we use all the remaining labels between $\frac{mn}{4}+1$ and $\frac{%
3mn}{4}$ for all the vertices in every odd column. Thus for $0\leq i\leq m-1$
and $j\in \{1,3,\ldots ,n-1\}$ let:

\begin{equation*}
\ell (v_{i,2j+1})=\left\{
\begin{array}{lcl}
\ell (v_{i,2j})+\frac{mn}{4}, & \text{if} & \ell (v_{i,2j})\leq \frac{mn}{2},
\\
\ell (v_{i,2j})-\frac{mn}{4}, & \text{if} & \ell (v_{i,2j})>\frac{mn}{2}.%
\end{array}%
\right.
\end{equation*}%
Obviously the labeling $\ell $ is a bijection from $V(C_{m}\times C_{n})$ to
$\{1,\ldots ,mn\}$. It is also straightforward to see that $k=2mn+2$ is the
magic constant. Hence $\ell $ is distance magic labeling.

However, $\ell$ is not balanced distance magic, as none of the
cycles $C_m, C_n$ is (see Theorem \ref{MRSC4}) and thus their
product cannot be balanced distance magic due to Theorem
\ref{str_pr}. ~\hfill \rule{0.1in}{0.1in}\medskip

The example of distance magic labeling of $C_{16}\times C_{16}$ is shown
below, where $v_{0,0}$ starts in lower left corner and the first index is
for the row and the second for the column:

\begin{equation*}
\begin{array}{cccccccccccccccc}
196 & 132 & 62 & 126 & 194 & 130 & 64 & 128 & 193 & 129 & 63 & 127 & 195 &
131 & 61 & 125 \\
228 & 164 & 30 & 94 & 226 & 162 & 32 & 96 & 225 & 161 & 31 & 95 & 227 & 163
& 29 & 93 \\
57 & 121 & 199 & 135 & 59 & 123 & 197 & 133 & 60 & 124 & 198 & 134 & 58 & 122
& 200 & 136 \\
25 & 89 & 231 & 167 & 27 & 91 & 229 & 165 & 28 & 92 & 230 & 166 & 26 & 90 &
232 & 168 \\
204 & 140 & 54 & 118 & 202 & 138 & 56 & 120 & 201 & 137 & 55 & 119 & 203 &
139 & 53 & 117 \\
236 & 172 & 22 & 86 & 234 & 170 & 24 & 88 & 233 & 169 & 23 & 87 & 235 & 171
& 21 & 85 \\
49 & 113 & 207 & 143 & 51 & 115 & 205 & 141 & 52 & 116 & 206 & 142 & 50 & 114
& 208 & 144 \\
17 & 81 & 239 & 175 & 19 & 83 & 237 & 173 & 20 & 84 & 238 & 174 & 18 & 82 &
240 & 176 \\
212 & 148 & 46 & 110 & 210 & 146 & 48 & 112 & 209 & 145 & 47 & 111 & 211 &
147 & 45 & 109 \\
244 & 180 & 14 & 78 & 242 & 178 & 16 & 80 & 241 & 177 & 15 & 79 & 243 & 179
& 13 & 77 \\
41 & 105 & 215 & 151 & 43 & 107 & 213 & 149 & 44 & 108 & 214 & 150 & 42 & 106
& 216 & 152 \\
9 & 73 & 247 & 183 & 11 & 75 & 245 & 181 & 12 & 76 & 246 & 182 & 10 & 74 &
248 & 184 \\
220 & 156 & 38 & 102 & 218 & 154 & 40 & 104 & 217 & 153 & 39 & 103 & 219 &
155 & 37 & 101 \\
252 & 188 & 6 & 70 & 250 & 186 & 8 & 72 & 249 & 185 & 7 & 71 & 251 & 187 & 5
& 69 \\
33 & 97 & 223 & 159 & 35 & 99 & 221 & 157 & 36 & 100 & 222 & 158 & 34 & 98 &
224 & 160 \\
1 & 65 & 255 & 191 & 3 & 67 & 253 & 189 & 4 & 68 & 254 & 190 & 2 & 66 & 256
& 192%
\end{array}%
\end{equation*}

Next theorem that completely describes distance magic graphs among
direct product of cycles follows immediately by
Theorems~\ref{notdm}, \ref{0mod4}, and \ref{str_pr}.

\begin{mytheorem}
A graph $C_{m}\times C_{n}$ is distance magic if and only if $n=4$ or $m=4$
or $m,n\equiv 0\imod 4$. It is
balanced distance magic if and only if $n=4$ or $m=4$.~\hfill \rule%
{0.1in}{0.1in}\label{main1}
\end{mytheorem}

\end{document}